\documentclass[12pt]{article}
\usepackage{amssymb}
\def\be{\begin{equation}}
\def\ee{\end{equation}}

\def\R{{\sf I\kern-.2em R}}
\def\N{{\sf I\kern-.2em N}}
\def\C{\kern.1em{\raise.47ex\hbox{$\scriptscriptstyle
$}}\kern-.40em{\sf C}}
\def\Z{{\sf Z\kern-.32em Z}}

\def\be{\begin{equation}}
\def\ee{\end{equation}}

\newtheorem{theorem}{\noindent Theorem}
\newtheorem{lemma}{\noindent Lemma}
\newtheorem{definition}{\noindent Definition}
\newtheorem{corollary}{\noindent Corollary}

\newtheorem{remark}{\noindent Remark}

\begin{document}

\begin{titlepage}

\begin{center}
{\LARGE \bf Distance matrices, random metrics and Urysohn space.}
\bigskip
\bigskip

{\large A.~M.~VERSHIK $^{1}$}
\bigskip

{\sl
$^{1}$ Steklov Institute of Mathematics at St.~Petersburg,\\
Fontanka 27,
119011 St.~Petersburg, Russia.\\
and
\sl Max-Planck-Institute fur Mathematik.

Partially supported by RFBR, grant 99-01-00098.}
\end{center}
\bigskip
\centerline{28.02.2002.}
\bigskip
\begin{center} ABSTRACT
\end{center}

   We introduce an universum of the Polish (=complete separable metric)
spaces - the convex cone of distance matrices and study its geometry.
It happened that the generic Polish spaces in this sense of this universum
is so called Urysohn spaces defined  by P.S.Urysohn in 20-th,
and generic metric triple (= metric space  with probability borel measure)
is also Urysohn space with non-degenerated measure.
We prove that the complete invariant of the metric space with measure up to measure preserving isometries is so called matrix distribution
 - a $S_{\infty}$-invarinat ergodic measure on the cone of distance matrices
 This defined an important new class of random matrices and family of random metrics on the naturals.

\begin{center} CONTENT
\end{center}

1.The cone of the distance matrices - the universum of the Polish spaces.

2.Geometry and topology of the cone of distance matrices.

3.Universal matrices and genericity of Urysohn space.

4.Metric triples and its invariant - matrix distribution.

5.Properties of matrix distributions, random metrics on the naturals.

Bibliography.

\end{titlepage}
\newpage

\section{The cone of distance matrices - the universum of the Polish spaces}

 Denote by $\cal R$ a set of all infinite real matrices
 $${\cal R} =\{\{r_{i,j}\}_{i,j=1}^{\infty}:  r_{i,i}=0, r_{i,j}\geq 0,
 r_{i,j}=r_{j,i},  r_{i,k}+r_{k,j}\geq r_{i,j},i,j,k=1 \dots\}$$

We will call such matrices (finite or infinite) {\it the distance matrices}.
Each such a matrix defines a semi-metric $\rho$ on an ordered countable set - for
definiteness - on the set of naturals $\bf N$: $\rho(i,j)= r_{i,j}$.
We allow zeros out from principle diagonal, so it is only semi-metric
in general. If matrix has no zeros out of principle diagonal we will call
it {\it true distance matrix}.

The set of all distance matrices is a weakly closed convex cone in the
real linear space $Mat_{\bf N}(\bf R)={\bf R^N}^2$ which equipped with ordinary weak topology.
We denote  this cone as ${\cal R}$ and will quote as the (topological)
space of distance matrices. Subset of true distance matrix is every dense open
strata in ${\cal R}$.

Suppose now that $(X, \rho)$ is a complete separable metric (=Polish) space
with metric $\rho$, and $\{x_i\}_{i=1}^\infty$ is an ordered dense countable set in it.
Define the matrix $r=\{r_{i,j}\} \in \cal R$ where
$r_{i,j}=\rho(x_j,x_j), i,j =1\dots$  is a distance matrix with positive elements
out from diagonal; we interpret it as a {\it metric on the set of naturals}.
Evidently this distance matrix includes all information about initial
space $(X,\rho)$, because  $(X,\rho)$ is canonical
completion of the set of naturals with that metric.

Any invariant metric property of the space (like compactness,
topological, homological properties etc,) could be expressed in terms of the
distance matrix for any dense countable subsets of that space. Some of them is easy to rewrite
(say, compactness) another more difficult (dimension).

General distance matrix (with possible zeros out of diagonal)
defines on the set naturals the structure of
{\it semi-metric space}. By completion of naturals in that case we mean the completion of
corresponding quotient metric space of the classes of points with zero distances.
Number of the classes we will call {\it geometric rank} of distance matrix or
of the corresponding semi-metric space. For example zero matrix is a distance matrix
of the naturals (or countable semi-metric space) with zero distance between each two points
and ``completion'' of it is one point metric space. So finite metric spaces also could
be considered in this setting.

 The space of distance matrices $\cal R$ is an
{\it universum of all separable complete metric spaces} with a fixed dense countable subsets:
we can look at $\cal R$ as a ``fiber bundle'', the base of which is
a set of all individual Polish spaces and the fiber over given space is
the set of all countable ordered dense subsets in this metric
space. Because of universality of Urysohn space $\cal U$(see below) the set of all
closed subsets of $\cal U$  could be considered as a base of that bundle.

Let $\xi$ is a partition of $\cal R$ on the classes of matrices which
produce the isometric completions of the set of naturals with that
metric as complete Polish spaces. The quotient space over  partition $\xi$
is the space of the {\it classes} of the Polish spaces up to isometry.
As was conjectured in  \cite{V} and proved in the paper \cite{CK01}
that partition (or equivalence relation) is not "smooth" so,
the quotient has no good borel or topological structure and
the problem of the classification of the Polish spaces up-to
isometry is wild (non smooth) problem. In the same time the restriction
of this problem to the case of compact Polish spaces is smooth (see \cite{G})
and a space of all classes up to isometry of compact metric spaces
has a natural topology.

The space $\cal R$  plays role of "tautology fibration" over
the space of classes of isometrcal
Polish spaces similar to ordinary topological constructions. We will
see that the problem of classification of the metric space with
borel probability measure is ``smooth'' problem and the complete invariant
is a measure on  $\cal R$ with a special properties.

In this paper we study geometry of the space $\cal R$ itself
and describe its generic properties Polish spaces in terms of that cone,
and prove that Urysohn space is generic.
Then we consider the metric triples (Gromov's mm-spaces), find
its invariant and prove that probability measures on the space $\cal R$.
which are concentrated on so called universal matrices are also generic,
so typical metric space with measure is also Urysohn space with measure.
Roughly speaking - random choice of the Polish space gives you that remarkable space.
More precise fromulation of this theorem will be done else where - this is the
analogues in another categories of Erdos-Renyi theorem (see (\cite{ER},{Ca})
about random graphs.
\newpage
\section{Geometry and topology of the cone ${\cal R}$}
\subsection{Convex structure}

Let us denote the finite dimensional cone of distance matrices of order $n$
as ${\cal R}_n = {\cal R} \cap Mat_n({\bf R})$. Cone ${\cal R}_n$ is polyhedral cone
inside the positive orthant in $ Mat_n({\bf R})\equiv {\bf R}^{n^2}$.
Denote as ${\bf M}^s_n({\bf R})\equiv {\bf M}^s_n$ the space
{\it symmetric matrices with zeros on the  principle diagonal}, the cone
 ${\cal R}_n$ is contained in this space: ${\cal R}_n \subset {\bf M}^s_n$
and the last space is evidently the linear hull of the cone:
$span({\cal R}_n)={\bf M}^s_n$, because the interior of ${\cal R}_n$
in that space is not empty. It is not so evident that $span({\cal R})={\bf M}^s_{\bf N}$,
where  ${\bf M}^s_{\bf N}$ is the space of all real infinite symmetric matrices with
zero principle diagonal (we do not need in weak closure of linear hull as it usually happened
in infinite dimensional cases), but we will not use this fact.

Each  matrix $r \in {\cal R}_n$ defines the (semi)metric
space $X_r$ on  the sets of $n$ points

Define projection $$p_{m,n}:{\bf M}^s_m \longrightarrow {\bf M}^s_n, m >n $$
which associates to the matrix $r$ of order $m$ its of NW-corner of order $n$.
The  cones ${\cal R}_n$ are preserved by $p_{n.m}: p_{m,n}({\cal R}_m)={\cal R}_n$
Projection $p_{n,m}$ are natural extends to the space of infinite
symmetric matrices with zero diagonal - $p: M^s_{\bf N}  \longrightarrow M^s_n({\bf R})$
and $p_n$ are also preserved the cones: $p_n({\cal R})={\cal R}_n$.
It is clear that ${\cal R}$ is inverse limit as topological space (in weak topology)
of  the finite dimensional cones of ${\cal R}_n$ under projections $\{p_n\}$.
We will omit the first index and denote $p_{N,n}=p_n$.

Let us consider a geometrical structure of ${\cal R}_n$ and ${\cal R}$.

It is easy to describe the extremal rays (in the sense of convex geometry)
of the convex polyhedral cone ${\cal R}_n, n=2, \dots, \infty$.
(for $n=1$ the cone ${\cal R}_1$ consist with one point - zero).

\begin{lemma}
{Each extremal ray in ${\cal R}_n, n>1$ is a ray of type
$\{\lambda \cdot l\}_{\lambda \geq 0}$,
where $\lambda$ runs over all real nonnegative numbers, and $l$ is a symmetric $0-1$-
distance matrix which corresponds to semi-metric space metric quotient of which has
just two points. The same is true for ${\cal R}$.

In another words all nonzero matrix $l$ which belongs to extremal ray
is a distance matrix of the finite or countable
\it semi-metric space which divides on two nonempty subset
the distance between two points from the same subset is zero and between two points
of the different subsets equal to one. The number of extremal rays of ${\cal R}_n, n>1$
is equal to $2^{n-1}-1$.}
\end{lemma}
\begin{proof}
{If the ray $\{{\lambda \cdot l}_{\lambda \geq 0}\}, l \in {\cal R}$ is extremal then the
non zero matrix $l$ must have at least one zero coordinates out from diagonal.
If we facotrize the space modulus zero distance
we obtain a new metric space which can not contain a non degenerated
triangles because existence of such triangles contradicts to extremality of the ray.
Consequently  the corresponding metric space has two points only.}
\end{proof}

The cones ${\cal R}_n$ as the topological spaces could be described
in the different ways. They have very interesting stratifications as semi-algebraic sets
and various systems of coordinates. In order to clarify topological structure of the cones
we will use the system of the coordinates which is well-organized as a sets of vectors
of entries of matrices, so called admissible vectors, which is very convenient
for the natural construction and studying of Urysohn space and measures on it
(see next paragraph).

\subsection{Admissible vectors and structure of the $\cal R$}

Suppose $r=\{r_{i,j}\}_1^n$ is a distance matrix of order $n$,
choose a vector $a \equiv \{a_i\}_{i=1}^n \in {\bf R}^n$
such that if we border a matrix $r$ with vector $a$ as the last
column and the last row then the {new matrix of order $n+1$ still will be
also a distance matrix. We will call such a vector {\it admissible vector}
for fixed distance matrix $r$ and denote the set of of all admissible vectors
for $r$ as $A(r)$. For given $a \in A(r)$ denote as $(r^a)$,
a distance matrix of order $n+1$ which is obtains from  matrix $r$
with adding vector $a \in A(r)$ as the last row and column.
It is clear that $p_n(r^a)=r$. The matrix $r^a$ has the form:

 $$\left( \begin{array}{ccccc}
 0    &r_{1,2}& \ldots&  r_{1,n}&a_1\\
 r_{1,2} &0& \ldots &r_{2,n} &a_2 \\
 \vdots &  \vdots&
\ddots & \vdots &\vdots\\
r_{1,n}&r_{2,n}& \ldots &0  &a_n\\
  a_1&  a_2  & \ldots & a_n& 0
\end{array}\right)$$

The (semi)metric space $X_{r^a}$ corresponding to matrix $r^a$ is extension of
of $X_r$ - we add one new point $x_{n+1}$ and $a_i, i=1 \dots n$ is a distance
between $x_{n+1}$ and $x_i$.
The admissibility  of $a$ is equivalent to the following set of inequalities :
vector $a=\{a_i\}_{i=1}^n$ must subtract to the series of
triangle inequalities
for all $i,j,k=1 \dots n$; (matrix $\{r_{i,j}\}_{i,j=1}^n$ is fixed):
\begin{equation}
\label{trivial} |a_i-a_j|\leq r_{i,j}\leq a_i+a_j
\end{equation}

So, for given distance matrix $r$ of order $n$ the set of admissible vectors
is $A(r)=\{\{a_i\}_{i=1}^n:|a_i-a_j|\leq r_{i,j}\leq a_i+a_j, i,j=1 \dots n\}$.
The set $A(r)$ is the intersection of cone ${\cal R}_{n+1}$ with affine subspace
which consists with matrices of order $n+1$ with matrix $r$ as a NW-corner of order $n$.
 It is clear from the linearity of inequalities that set $A(r)$ is an
unbounded closed convex polytope in ${\bf R}^n$.
If $r_{i,j} \equiv 0, i,j=1 \dots n\geq 1$,
then $A(r)$ is diagonal: $A(r)=\Delta_n \equiv \{\lambda,
\dots (n) \dots \lambda\}_{\lambda \geq 0}\subset {\bf R}^n_+$.
Let us describe its structure more carefully.
\begin{lemma}
{For each true distance matrix $r$ of order $n$
the set of admissible vectors $A(r)$
is a closed convex polyhedron in orthant ${\bf R}_+^n$, namely this is a Minkowski sum:
$$A(r)= M_r +\Delta_n,$$
where  $\Delta_n \equiv \{\lambda,
\dots (n) \dots \lambda\}_{\lambda \geq 0}$ is a half-line -- positive diagonal in the
space ${\bf R}_+^n$,
 and $M_r=$ conv(ext $A(r))$  is a compact convex polytope of the dimension
$n$. This polytope is a convex hull of extremal points of the set $A(r)$.

The set $A(r)$ is homeomorphic to the product of the simplex and half-line,
or, equivalently  to the half-space  $R_+^n=\{(b_1, \dots b_n): b_n \geq 0\}$)
(this homeomorphism could be chosen as piecewise linear but not canonically).
 If $r$ is not true distance matrix of order $n$ and has geometric rank $m<n$
then $A(r)$ is homeomorphic to the product of $(m-1)$-dimensional simplex and half-line.}
\end{lemma}
\begin{proof}
{The set $A(r) \subset {\bf R}^n$ is intersection of finitely many of the closed
subspaces, evidently it does not contain straight lines, so, because of
general theorem of convex geometry $A(r)$ is a sum of the convex closed polytope
and some cone with the vertex at origin.
That convex polytope is the convex hull of the extremal points of convex set $A(r)$.
But this cone must be one dimensional, namely - diagonal in ${\bf R}^n$ because
if it contains any half-line differ from diagonal then
the triangle inequality (left side) will violate. Dimension of  $A(r)$
is a number which depends on matrix $r$ and could be less $n$ for degenerated $r$,
evidently dimension of $M_r$ is equal to dimension of $\dim A(r)$ or to $\dim A(r)-1$.
The assertion about topological structure
of $A(r)$ follows from what was claimed above.}
\end{proof}

The convex structure of polytopes $M_r, A(r)$ is very interesting and it seems never
had been studied. For dimensions more than 3
combinatorial type of the polytope $M_r$ hardly depends on $r$. For our purpose here
it is important only to establish {\it topological} isomorphisms of the special type
with half-spaces for true distance matrices of the given order.

In dimension three the {\it combinatorial type} of polytopes $M_r$, and consequently
combinatorial structure of the sets $A(r)$ is the same for all true distances matrices $r$.

{\bf Example}
For $n=3$ the description of the set $A(r)$ and its of extremal points is the following.
Let $r^a$ is the matrix:

 $$\left( \begin{array}{cccc}
 0    &r_{1,2}&   r_{1,3}&a_1\\
 r_{1,2} &0&r_{2,3} &a_2 \\
 r_{1,3}&r_{2,3}&0  &a_3\\
  a_1&  a_2  &  a_3 &  0
\end{array}\right)$$
There are seven extremal points $(a_1,a_2,a_3)$ of $A(r)$ :
 the first one is a vertex which is the closest to origin:
$(a_1,a_2,a_3)=(\frac{1}{2}(r_{1,2}+r_{1,3}-r_{2,3}, \frac{1}{2}(r_{1,2}-r_{1,3}+r_{2,3}),
\frac{1}{2}(-r_{1,2}+r_{1,3}+r_{2,3}))$;

three non degenerated extremal points:
$(\frac{1}{2}(r_{1,2}+r_{1,3}+r_{2,3}), \frac{1}{2}(-r_{1,2}+r_{1,3}+r_{2,3}),
\frac{1}{2}(r_{2,3}+r_{1,3}-r_{2,3}))$ and the rest two are cyclic permutations of that
in the natural sense.

and  three degenerated extremal points
$(r_{1,2},0,r_{2,3}), (0, r_{1,2}, r_{2,3}), (r_{1,3}, r_{2,3}, 0)$
which defined  the metric spaces with
two distinguish points and with the third point which coincide with one of those.

If $r_{1,2}=r_{1,3}=r_{2,3}=1$
then those seven points are as follow

$(1/2,1/2,1/2), (1,0,1),(0,1,1),(1,1,0),
(3/2,1/2,1/2),(1/2,3/2,1/2),(1/2,1/2,3/2)$.

It makes sense to mention that all non-degenerated extremal points defines
the finite metric spaces which can not be isometrically embedded to Euclidean space.
 The cone ${\cal R}_3$ is simplicial cone of dimension 3, and for all matrices
$r \in {\cal R}_3$ which does not belong to the boundary the set of admissible
vectors $A(r)$  (extremal points of which we had defined above) have the same
convex combinatorial type although they are no affine isomorphism between them.

\subsection{Projections and isomorphisms}

Suppose we have distance matrix $r$ of order $N$ and its NW-corner of order $n<N$
- $p_n(r)$. Then we can define neglecting projection
$\chi^r_n$ of $A(r)$ to $A(p_n(r)$:
$\chi^r_n:(b_1, \dots b_n, b_{n+1},\dots b_N) \mapsto (b_1 \dots b_n)$.
(We omit index $N$ in the denotation of $\chi^r_n$)
The next simple lemmas play very important for our construction.
\begin{lemma}

{Let $r\in {\cal R}_n$ is a distance matrix of order $n$ and two vectors
$a=(a_1, \dots a_n) \in A(r), b=(b_1, \dots b_n) \in A(r)$ there exist
real nonnegative number $h \in {\bf R}$ such that vector
${\bar b}=(b_1,\dots b_n,h) \in A(r^a)$.}
\end{lemma}

The claim of this lemma is equivalent to the assertion that for each
$r$ the projection $\chi^r_n$ defined above is epimorphism of $A(r^a$ to $A(r)$:

\begin{corollary}
{For each $r \in {\cal R}_n$ and $a \in A(r)$ the map:
$\chi^r_{n+1,n}: (b_1, \dots b_n, b_{n+1}) \mapsto (b_1 \dots b_n)$ of  $A(r^a) \to A(r)$
is epimorphism of $A(r^a)$ on $A(r)$.(by definition $p_{n+1,n}(r^a)=r$)}
\end{corollary}

\begin{proof}
{The assertion of lemma as we will see, is a simple geometrical observation: suppose we have
 two finite metric space
$X=\{x_1, \dots x_{n-1}, x_n\}$ with metric $\rho_1$
and $Y=\{y_1, \dots y_{n-1},y_n \}$ with metric $\rho_2$. Suppose the subspaces
of the first $n-1$ points $\{x_1, \dots x_{n-1}\}$ and $\{y_1, \dots y_{n-1}\}$
are isometric -- $\rho_1(x_i,x_j)=\rho_2(y_i,y_j), i,j=1, \dots n-1$,

then there exists the third space $Z=\{z_1,\dots z_{n-1},z_n,z_{n+1}\}$ with metric
$\rho$ and two isometries  $I_1,I_2$ of both spaces $X$ and $Y$ to the space
$Z, I_1(x_i)=z_i,I_2(y_i)=z_i, i=1, \dots n-1, I_1(x_n)=z_n,
I_2(y_n)=z_{n+1}$.
In order to prove existence of $Z$ we need to prove that it is possible to
define only nonnegative number $h$ which will be
the distance $\rho(z_n,z_{n+1})=h$ between $z_n$ and $z_{n+1}$ (images of $x_n$ and $y_n$
in $Z$ correspondingly)
such that all triangle inequalities
took place in the space $Z$.
The existence of $h$ follows from the inequalities:
$$\rho_1(x_i,x_n)-\rho_2(y_i,y_n) \leq \rho_1(x_i,x_n)+\rho_1(x_i,x_j)-\rho_2(y_i,y_n)=$$
$$=\rho_1(x_i,x_n)+\rho_2(y_i,y_j)-\rho_2(y_i,y_n)\leq \rho_1(x_i,x_n)+\rho_2(y_i,y_n)$$
for all $i,j=1, \dots n-1$.

 Consequently $$\max_i |\rho_1(x_i,x_n)-\rho_2(y_i,y_n)| \equiv m
\leq  M \equiv  \min_j (\rho_1(x_i,x_n)+\rho_2(y_i,y_n)).$$
So, a number $h$ could be chosen as an arbitrary number from the nonempty closed interval
$[m,M]$ and we define $ \rho(z_n,z_{n+1}) \equiv h$;
it follows from the definitions that all triangle inequalities are satisfied.
Now suppose we have a distance matrix $r$ of order $n-1$ and admissible vector $a \in A(r)$,
so we have metric space $\{x_1, \dots x_{n-1}, x_n\}$ (first $n-1$ points corresponds
to matrix $r$ and all space - to extension matrix $r^a$. Now suppose we choose
another admissible vector $b \in A(r)$, distance matrix $r^b$ defined space
$\{y_1, \dots y_{n-1},y_n\}$ where subset of first $n-1$ points is isometric (the same)
as space $\{x_1, \dots x_{n-1}\}$. As we proved we can define space $Z$ whose
distance metric $\bar r$ of order $n+1$ gives has needed property.}
\end{proof}

Now we can formulate the general assertion about projections $\chi^r$.

\begin{lemma}
{For each naturals $N,n$ and $r \in {\cal R}_N$  the map $\chi^r_n$ is epimorphism
of $A(r))$ onto $A(p_n(r))$. In another words for each $a=(a_1, \dots a_n) \in A(p_n(r))$
there exist vector $(b_{n+1},\dots b_N)$ such that
 $b=(a_1,\dots a_n,  b_{n+1}, \dots b_N) \in A(r)$.}
\end{lemma}
\begin{proof}
{Previous proof shows how to define the first number $b_{n+1}$.
But the projection $\chi^r_n$ as a map from  $A(r), r\in {\cal R}_N$ to $A(p_n(r)$
is the product of projections $\chi^r_n\cdots \chi^r_{N-1}$
Previous lemma shows that all those factors are epimorphisms.}
\end{proof}

It is useful to  consider each infinite distance matrix
$r \equiv \{r_{i,j}\} \in {\cal R}$  as
as a sequences of the admissible vectors of the increasing lengths
\begin{equation}
r(1)=\{r_{1,2}\},
r(2)=\{r_{1,3},r_{2,3}\}, \dots r(k)=\{r_{1,k+1},r_{2,k+1}, \dots r_{k,k+1}\}\dots,
k=1,2 \dots,
\end{equation}
with conditions
$r(k) \in A(p_k(r))$, (remember that $p_k(r)$ is NW-projection of matrix $r$ on the space
${\bf M}^s_k$ which was defined above) -  each vector $r(k)$ is admissible
for the {\it previous distance matrix}. This is a realization of the infinite distance matrix
using the projections of the cones. We can consider the following sequences
of the cones and maps:

\begin{equation}
\label{trivial}0={\cal R}_1\stackrel{p_2}{\longleftarrow}
{\cal R}_2={\bf R}_+\stackrel{p_3}{\longleftarrow}{\cal R}_3
{\longleftarrow} \dots
{\longleftarrow}{\cal R}_{n-1}\stackrel{p_n}{\longleftarrow}{\cal R}_n{\longleftarrow} \dots
\end{equation}
where projections $p_n$ now is restriction of the projection $p_n$ (defined above)
onto the cone ${\cal R}_n$. A preimage of the point $r \in {\cal R}_{n-1}$ (fiber over $r$)
 is the set $A(r)$ which structure had been described in the lemma 2. The fibers are not
even homeomorphic to each other for the various $r$ (even dimensions could be different),
so the sequences (3)
are not the sequences of fiber bundles in usual sense. But it defines a very interesting
stratifications of all cones, and the corresponding complex structure on $\cal R$.
More careful studying of the stratification of the cones ${\cal R}_n$ and  ${\cal R}$
must take in account the combinatorial and semi-algebraic structures of its.
We postpone the discussion on this interesting subject.

 For our goals in this paper
it is enough to consider the open cell (main strata) of each cones which is simply
{\it the set of distance matrices for which  all triangle inequalities
are strict inequalities}.
Denote this open part (interior) of the cone ${\cal R}_n$ as ${\cal R}_n^0$.
The cone ${\cal R}^0$ is a set of  all distance matrices for which any triangle inequalities
are strict, it is every dense subset of $\cal R$
Let $r \in  {\cal R}_n^0$, denote as ${\dot A}(r)$
the interior of $A(r)$ (of those admissible vectors for which again
all triangle inequalities are strict).
So the cones ${\cal R}^0$ equipped with the structure of inverse limit.
For each $n$ projection  ${\cal R}_{n-1}^0\stackrel{p_n}{\longleftarrow}{\cal R}_n^0$
this is a trivial fibration because now all fibers are homeomorphic
to each other and to open half-space, and the base is homeomorphic to affine space,
because this is an open nonempty convex set. We can to refine this picture with
the trivialization of this fibration because now all fibers for given $n$ are
open cells of the same dimension.We obtain

\begin{lemma}
{1.The structure of the cone ${\cal R}^0$ is described with the sequence of the maps:
\begin{equation}
\label{trivial}0={\cal R}_2^0={\bf R}_+\stackrel{p_3}{\longleftarrow}{\cal R}_3^0
{\longleftarrow} \dots
{\longleftarrow}{\cal R}_{n-1}^0\stackrel{p_n}{\longleftarrow}{\cal R}_n^0{\longleftarrow}
 \dots{\longleftarrow} {\cal R}^0
\end{equation}

For each $n$ the map
 ${\cal R}_{n-1}^0\stackrel{p_n}{\longleftarrow}{\cal R}_n^0$
defines a trivial fiber bundle and natural decomposition:
${\cal R}_n^0 \simeq{\cal R}^0_{n-1} \times  {\bf \dot R}_+^{n-1}$
where ${\bf \dot R}_+^{n-1}$ is open half-space, $\pi_n$
is the projection on the first summand in that decomposition;
and the second factor is homeomorphic to ${\dot A}(r)$ for all $r$.
So we have a realization of open strata:
$${\cal R}^0 \simeq \prod_{k=1}^{\infty} {\bf \dot R}_+.$$}
\end{lemma}

\begin{remark}
{The topology which is induced on ${\cal R}^0$ from ${\cal R}$ does not coincide
with topology of inverse limit in (2), but we will not use it.
The borel structure defined by the family of cylindric sets
on that direct product coincides with borel
structure which is induced on ${\cal R}^0$ from  $\cal R$.
There is no canonical isomorphism in the theorem: there are many ways to identify
the fibers with  affine open half-spaces (or simply with open discs) according to lemma 2.}
\end{remark}

\newpage

\section{Universality and Urysohn space}
\subsection{Universal distance matrices}
The following definition plays a crucial role.
\begin{definition}
{The infinite true distance matrix $r=\{r_{i,j}\}_{i,j=1}^{\infty} \in {\cal R}$
called as universal matrix if the following condition is satisfied:

for each $\epsilon >0$, each $n \in \bf N$ and for each
vector $a=\{a_i\}_{i=1}^n \in A(p_n(r))$
there exists such $m \in \bf N$  that
$\max_{i=1 \dots n}|r_{i,m}-a_i| < \epsilon.$

In another words: for each $n \in \bf N$ the set of vectors
$\{\{r_{i,j}\}_{i=1}^n\}_{j=n+1}^{\infty}$ is every dense in the set
of admissible vectors  $A(p_n(r))$.}
\end{definition}

Let us denote the set of universal distance  matrices as $\cal M$
\begin{theorem}
{The set of universal matrix in $\cal R$ is nonempty and moreover - is an every dense
$G_{\delta}$-set in weak topology of the cone $\cal R$.}
\end{theorem}

\begin{proof}
{We will use representation of the lemmas of previous section for the construction
of at least one universal true distance matrix in the cone ${\cal R}_0$.

Let us fixed sequence  $\{m_n\}_{n=1}^\infty$ of naturals in which each
natural number occurs infinitely many times and with property: for each $n,  m_n \leq n, m_1=1$.
For each finite distance matrix $r \in {\cal R}_n$ let us choose an ordered countable dense
subset
$\Gamma_r=\{\gamma^r_k\}_{k=1}^\infty \subset \dot A(r){\subset \bf R}^n$ and
choose any metric in  $\dot A(r)$, say Euclidean norm from $\bf R^n$.

The first step consists with the choice of positive real number
$\gamma_1^1 \in \Gamma_1 \subset {\dot A}(0)={\bf \dot R}_+^1$
so we define a distance matrix $r$ of order $2$ with element $r_{1,2}=\gamma_1^1$.

Our method of the construction of the  universal matrix $r$ used its representation
as sequence of the admissible vectors $\{r(1),r(2),\dots\}$  of increasing lengths (see (2)),
- the index in the brackets is a dimension of the vector, -the conditions
on the vectors are as follow: $r(k) \in A(p_k (r_{k+1}))$. The
sequence of the corresponding matrices $r_n, n=1 \dots $
are stabilized to the infinite matrix $r$.
Suppose  after  $n-1$-th step we obtain finite matrix $r_{n-1}$,
now we choose a new admissible vector $r(n) \in {\dot A}(r_{n-1})$. The choice of
this vector (denote it as $a$) is defined by condition: the
distance in $A(r_{m_{n}})$ (in norm) between projection $\chi^r_{m_n}(a)$
of the vector $a$ onto subspace of admissible vectors $A(r_{m_n})$
and point $\gamma^{m_n}_s \in \Gamma_{r_{m_n}} \subset A(r_{m_n}) $
must be less than $2^{-n}$, here $s =|i: m_i=m_n, 1\leq i \leq n|+1$:
$$\|\chi^r_{m_n}(a) - \gamma^{r_{m_n}}_s\|<2^{-n}.$$
 Remain that projection $\chi^r_{m_n}$ is epimorphism from $A(r)$ to $A(p_{m_n}(r))$,
(lemma 4), so a vector $a \in A(r_n)$ with this properties does exist.
Number $s$ is nothing more than the numbers of the points of $\Gamma_{r_{m_n}}$
 which occur on the previous steps of the construction.
After countable many of steps we obtain the infinite distance
matrix $r$.

 Universality of $r$ is evident, because for each $n$ projection $\chi^r_n$
 of vectors $r(k), k=n+1 \dots$  is a dense set in $A(r_n)$ by construction.
 This proves the existence of universal matrix.

From the previous proof it is easy to extract by induction over $n$
the following very important property of any universal matrix $r$.
Suppose $r=\{r_{i,j}\} \in \cal R$ is an arbitrary universal distance matrix,
and $q$ is a finite distance matrix of order $n \geq 2$ and $\epsilon >0$ then there
exists a submatrix $\{r_{i_s,i_t}\}_{s,t=1}^n$ of matrix $r$, such that
$\|q-{\bar q}\|< \epsilon$
(here $\|.\|$ -is an arbitrary norm in the space of matrices of order $n$).
 Indeed, because $r^1=\{r_{1,1\}}=0$ then $A(r^1)={\bf R}_+$ and by universality
$\{r_{1,n}\}_n^\infty$ must be dense in ${\bf R}_+ $, so we can choose
$m_1$ such that $|r_{1,m_1}-q_{1,2}|<\epsilon$, then because of density
of the columns of length 2 we can choose $m_2$ such that
$\|(r_{1,m_2},r_{2,m_2})-(q_{1,3},q_{2,3})\| < \epsilon$ etc.

We can say that for each $n$ the closure of the set of all submatrices
 $q= \{r_{i_s,i_t}\}_{s,t=1}^n$ of order $n$ over all choice of
n-plies $i_1<i_2< \dots <i_n$ of universal distance matrix
$r=\{r_{i,j}\}$ coincides with the set ${\cal R}_n$ of
all distance matrices of order $n$. Let us call this property a {\it weak universality}.

Now remark that the universality of the matrix is preserved under the action
of any finite permutations which simultaneously permute rows and columns,
also universality is preserved under the NW-shift which cancels the first row
and first column of the matrices.
Finally the set of universal matrices $\cal M$ is invariant under the changing
of the finite part of the matrix. (or - universality is {\it stable property}).
Consequently, $\cal M$ contains with
the given matrix all its permutations and shifts. But because of the weak
universality, the orbit of $r$ under group of permutations $S_{\bf N}$
is every dense in $\cal R$ in weak topology.

Finally, the formula which follows directly from the
definition of universality shows us immediately that
the set of all universal matrices $\cal M$
is a $G_{\delta}$-set:

$${\cal M} =\cap_{k \in {\bf N}}\cap_{n \in {\bf N}}
\cap_{a \in A(r^n)} \cup_{m \in {\bf N}, m>n}\{r\in {\cal R}:\max_{i=1,\dots n}
|r_{i,m}-a_i|< \frac{1}{k} \}.$$}
\end{proof}
\begin{remark}
 The property of matrix to be universal is much more stronger than the property
of weak universality  - it is easy to give an example of non universal
but weak universal matrix.
\end{remark}

\subsection{Urysohn space}
Now we introduce a remarkable Urysohn space.
In one of his last papers \cite{U} Urysohn gave a concrete construction of the
universal Polish space which we will call "Urysohn space" and denote as $\cal U$.
There are no the notion of universality in \cite{U} because Urysohn did not consider
infinite matrices at all but he actually have proved
several theorems which we summarize as following theorem:

\begin{theorem}(Urysohn \cite{U})
{There exist the Polish space with the properties:

1)(Universality) For each Polish space $X$ there exists the isometric embedding to $\cal U$;

2)(Homogeneity) For each two isometric finite subsets $X=(x_1 \dots x_m)$
and $Y=(y_1 \dots y_m)$ of $\cal U$ there exists isometry $I$ of whole space $\cal U$
which brings $X$ to $Y$;

3)(Uniqueness) Polish spaces with the previous properties is unique up to
isometry}
\end{theorem}

In our terminology Urysohn actually have proved
\begin{theorem}
{Let $r$ is an arbitrary infinite universal distance matrix, consider the set of
naturals $\bf N$ with metric $r$. Then the completion of the space
$({\bf N},r)$ with respect to metric $r$ is Urysohn space}
\end{theorem}

The proof of all three items of the theorem  follows directly and easily from
the  universality of the matrix $r$ and does not use even the existence of universal matrices.
That was a content of three main theorems in \cite{U}.
We do not repeated here very natural Urysohn's arguments and restrict ourself
only with the hints.
Universality of the matrix $r$ allows to construct a family
of fundamental sequences such their limits (in completion) gives the
finite set of points  with the given a prori distance matrix.
This gives universality (item 1). If we have two given finite isometric sets of
the points of completion $(x_1, \dots x_n)$ and  $(y_1, \dots y_n)$
we can include each of them to two countable dense subsets which are still
isometric - this gives homogeneity (item 2). Uniqueness follows from the fact
that in each space with the properties 1),2) it is possible to find
countable dense subset with given universal matrix.

\begin{remark}
{1.The inverse  assertion is trivially true: the distance matrix of any
countable ordered dense set in Urysohn space is an universal matrix.
In another words, universality of the true distance matrix $r$
is necessary and sufficient condition for the metric on the
naturals $({\bf N},r)$ to have as a completion Urysohn space.

2.Remark also, that if we use  in our construction admissible vectors
whose coordinates less or equal to 1 then we obtain
{\it universal metric space of diameter 1} with the same
property of universality for Polish space with diameter $\leq 1$,
this universal space was also mentioned in \cite{U,G}.}
\end{remark}

 Urysohn's paper begun with a very concrete example of the countable
space with rational distance matrix (indeed that was an universal
space over rationales - $\bf Q$).
More simple iterative method of the construction of the universal space
was suggested  in 80-th by Katetov, the similar idea was
used later by Gromov  in his book \cite{G}. Our method of construction
of the universal matrix differs from all those, - it used geometry
of the cone $\cal R$ and as we will see in the next section
allows interpret all the construction and  Urysohn space
in probabilistic. It also allows to study genericity of various
property ties of Polish spaces.

 As a result of the theorem 1 we obtain the remarkable fact that
"typical" (generic) distance  matrix is universal matrix, and {\it typical
(generic) in our unversum $\cal R$ Polish metric space is Urysohn space $\cal U$}!
In the next section we will prove the analog of this fact for metric spaces with measure.

\newpage
\section{Metric triples and its invariant - matrix distribution.}

Suppose $(X,\rho,\mu)$ Polish space with metric $\rho$ and with
borel probability measure $\mu$. We will call {\it metric triple} (In \cite{G} the author used
term ``mm-space''). Two triples  $(X_1,\rho_1, \mu_1)$ and
$(X_2, \rho_2, \mu_2)$ are {\it equivalent} if there exists isometry
$V$ which preserve the measure $$\rho_2(Vx,Vy)=\rho_2(x,y),  V \mu_1=\mu_2$$.
We have  mentioned before that the classification of the Polish space is non smooth
 problem. But the classification of metric triple has a good answer.

For triple  $T=(X,\rho, \mu)$ define the infinite product and Bernoulli measure
$(X^{\bf N},\mu^{\bf N}$
and define the map $F: X^{\bf N} \to {\cal R}$ as follow:

$$ F_T(\{x_i \}_{i=1}^\infty) =\{\rho(x_i,x_j)\}_{i,j=1}^\infty \in {\cal R}$$
The $F_T$-image of the measure $\mu^{\bf N}$ which we denoted as $D_T$ will called
{\it matrix distribution of the triple $T$} :$F_T \mu^{\infty} \equiv D_T$.

\begin{lemma}
{Measure $D_T$ is  borel probability measure on  $\cal R$
which is invariant and ergodic with respect to group of simultaneous permutations
of rows and columns and to simultaneous shift in vertical and horizontal directions
(NW-shift).}
\end{lemma}
\begin{proof}
{All facts follow from the same properties of measure $\mu^{\infty}$
and because map $F_T$ commutes with action of the shifts and permutations.}
\end{proof}

Let us call a measure on the metric space non-degenerated if there are no nonempty
open set of zero measure.
\begin{theorem}
{Two metric triples $T_1=(X_1,\rho_1, \mu_1)$ and $(X_2,\rho_2,\mu_2)$ with
non-degenerated measures  are equivalent iff its matrix distributions are coincided:
$D_{T_1}=D_{T_2}$ as the measures on the cone $\cal R$.}
\end{theorem}

\begin{proof}
{The necessity of the coincidence of the matrix distribution is evident:
if there exists the isometry $V$ between $T_1$ and $T_2$ which preserves measures then
the infinite power $V^{\infty}$ preserve the measures: $V^{\infty}(\mu_1^{\infty})=
\mu_2^{\infty}$ and because of equality $F_{T_1}X_1^{\infty}=F_{T_2}(V^{\infty}X_1^{\infty}$
the image of these measure are the same $D_{T_2}=D_{T_1}$.
Suppose now that  $D_{T_2}=D_{T_1}=D$. Then $D$-almost all distance matrices
$r$ are the images under the maps $F_{T_1}$ and $F_{T_2}$ say
$r_{i,j}=\rho_1(x_i,x_j)=\rho_2(y_i,y_j)$ but this means that the identification
of $x_i \in X_1$ and  $y_i\in X_2$ for all $i$ is isometry $V$ between these countable
sets $V(x_i)=y_i$. Now the main argument: by ergodic (with respect to NW-shift) theorem
$\mu_1$-almost all sequences of $\{x_i\}$ and $\mu_2$-almost all sequences
$\{y_i\}$ are uniformly distributed on $X_1$ and $X_2$ correspondingly.
This means that the $\mu_1$ measure of each ball
$B^l(x_i) \equiv \{z\in X_1:\rho_1(x_i,z)<l\}$
equal to $$\mu_1(B^l(x_i))
=\lim_{n \to \infty}\frac{1}{n} \sum_{k=1}^n 1_{\rho_1(x_k,x_i)\geq l} $$
which is equal because of isometry $V$ ($r_{i,j}=\rho_1(x_i,x_j)=\rho_2(y_i,y_j)$
- see above) to the $\mu_2$ measure of the ball:
$B^l(y_i) \equiv \{u\in X_2:\rho_2(y_i,u)<l\}$
 $$\mu_2(B^l(y_i))
=\lim_{n \to \infty}\frac{1}{n} \sum_{k=1}^n 1_{\rho_2(y_k,y_i)\geq l}= \mu_1(B^l(x_i)) $$
But both measures are non-degenerated, consequently each sequences $\{x_i\}$
and $\{y_i\}$ are every dense in its own space and because
both measures are borel it is enough to conclude their coincidence
if we established that the measures of the all such balls are the same.}
\end{proof}

\begin{corollary}
{Matrix distribution is complete invariant of the equivalence of metric triples
with non-degenerated measures.}
\end{corollary}

 Firstly this theorem as ``Reconstruction theorem'' in another formulation
had been proved in the book \cite{G} pp.117-123 by Gromov. He formulated its
in the terms of finite dimensional distributions of what we called
matrix distribution and proved using analytical method. He asked me
in 1997 about this theorem and I suggested  the proof which
is written here (see also in \cite{V}) and which he had quoted (pp.122-123) in the book.
Gromov has invited the readers to compare two proofs, one of which is
rather analytical (Weierstrass approximations) and another (above)
in fact uses only the reference on the
ergodic theorem. One of the explanations of this
difference is the same as in all applications of ergodic theorem -
it helps us to substitute the methods of space approximation onto operations
with infinite (limit) orbits. In our case the consideration of infinite
matrices and cone $\cal R$ with invariant measures gives a possibility
to reduce the problem to the investigation of ergodic action of infinite groups.
\newpage
\section{Properties of matrix distributions and random metrics on the naturals}

A matrix distribution of the metric triple $T=(X,\rho, \mu)$
is by definition a measure  $D_T$ on the cone $\cal R$.
It means that we can consider it as random (semi)metric
on the naturals. Namely, we choose independently the points $\{x_i\}$ of
metric space $X$ with measure $\mu$ as a distribution and
and look at their distance matrix $\{r_{i,j}\}=\{\rho(x_i,x_j)\}$.

Each distance matrix define the metric on $\bf N$ and consequently
the Polish space which is completion of naturals with that metric.
But if we have some probability measure on $\cal R$ then we need to
put some condition on this measure in order it corresponds to
some metric triple, or in another words - in order almost with respect
to that measure distance matrices define the same
(isometric) completion of naturals.

As we mentioned all measure $D_T$ must be invariant and ergodic with respect to
action of infinite symmetric group and to NW-shift. But this is not
 enough and not all restriction to which it satisfies.
We mention necessary and sufficient conditions (see \cite{V})
for that.

\begin{theorem}
{Let $D$  - a probability measure on the cone $\cal R$,which is invariant
and ergodic with respect to action of infinite symmetric group
(=group of all finite permutations of the naturals). The following
condition is necessary and sufficient in order $D$ is matrix distribution
for some metric triple $T=(X,\rho, \mu)$ or $D=D_T$:

for each $\epsilon>0$ there exists integer $N$ such that

$$D\{r=\{r_{i,j}\} \in {\cal R}: \lim_{n \to \infty} \frac{|\{j:1\leq\ j \leq n,
\min_{1 \leq i \leq N}
 r_{i,j} <\epsilon \}|}{n}> 1- \epsilon \}>1-\epsilon$$

    The following more stronger condition is equivalent to the compactness of
metric space: if for each $\epsilon >0$ there exists integer $N$ such that
$$D\{r=\{r_{i,j}\} \in {\cal R}: \mbox{for all} j>N,
\min_{1\leq i \leq N} r_{i,j} <\epsilon\}>1-\epsilon,$$
then $D=D_T$ where $T=(X, \rho, \mu)$ and $(X, \rho)$ is compact.}
\end{theorem}

The second condition is simply retelling of the existence of $\epsilon$-net for
each positive $\epsilon$ in the metric compactum;
the first one expresses well-known property of metric
space with measure: the set of full measure is sigma-compact,
or for each $\epsilon$ there exist a compact of measure $>1-\epsilon$.
The subset in the outer brackets in the first condition
is rather complicate borel set.

\begin{remark}
{As it follows from the general theorem about classification of the measurable functions
(see \cite{V2}), description of the  $S_{\infty}$-invariant measures
and previous theorem any measure which is concentrated on $\cal R$,
is $S_{\infty}$-invariant, ergodic, and simple - \cite{V2} -
is already a matrix distribution for some
metric on metric space with measure and consequently satisfied to the above condition.}
\end{remark}
The detales will be considered elsewhere.

Now let us consider the arbitrary measures on the cone $\cal R$,
or - an arbitrary random metrics on the naturals and choose some denotations.

{\bf Denotations.}
Denote  $\cal V$ the set of all probability borel measures \footnote{we use later on
term ``measure'' in the meaning - ``borel probability measure''} on the cone $\cal R$
and equip it with weak topology,- this is a topology of inverse
limit of the sets of probability measures on the finite dimensional cones ${\cal R}_n$
with its usual weak topology. The convergence in this topology is
convergence on the cylindric sets with open bases.

 Let $\cal D$  be  the subset in $\cal V$ of the  measures each of which is
a matrix distributions for some metric triple and equip it with induced topology
from $\cal V$. This set was described above in the theorem 5 in a constructive way
and, as claimed in the last remark, coinsides with $S_{\infty}$-invariant and ergodic
simple measures from $\cal V$.

The set $\cal P$ is a subset of  $\cal V$ of measures each of which concentrated
on the set of universal distance matrices:
$\nu \in \cal P$ iff $\nu({\cal M})=1$, and finally

${\cal Q} \subset \cal D$ consists with the measures which corresponds to the metric triples
$T=(X, \rho, \mu)$ in which $(X,\rho)$ is Urysohn space.
We have imbeddings
$$ {\cal P} \subset {\cal V} \supset {\cal D} \supset {\cal Q} \subset {\cal P};
      {\cal Q}= {\cal P}\cap{\cal D}  $$

Remark that the set of non-degenerated measures is of course every dense $G_\delta$
set in $\cal V$

\begin{theorem}
{1.The subset ${\cal P} \subset \cal V$
is every dense $G_\delta$ in $\cal V$. This means that for generic measures $\nu$
on $\cal R$
$\nu$-almost all distance matrices $r \in \cal R$ for $\nu \in \cal P$
are universal and consequently $\nu$-almost all distance matrix  $r$ defines the
metric on the naturals such that completion of naturals with respect to this
metric is the Urysohn space.

2.The subset ${\cal Q} \subset \cal D$
is every dense $G_\delta$-set in $\cal D$. This means that generic
metric triple $T=(X, \rho, \mu)$ has Urysohn space as the space $(X,\rho)$.}
\end{theorem}

\begin{proof}
{The first claim follows from the
 theorem 1 which claimed in particular
that the set of distance matrices which are universal is a set of type $G_\delta$
in $\cal R$, and from a general fact that the set all measures
on separable metrizable space such that the given every dense $G_\delta$-set
(in our case - $\cal M$)
has measure 1, is in its turn every dense $G_\delta$-set in the space
of all measures in weak topology.
The second claim follows from the fact that intersection
of $G_\delta$-set with any subspace in Polish space is also $G_\delta$-set in
induced topology.}
\end{proof}

Now we can give a probabilistic interpretation of our construction of
universal matrix and Urysohn space. Namely, we will give a family of measures from
the set $\cal P$.

{\bf Example.}

  Let $\gamma$ is a continuous measure
on open half-line ${\bf \dot R}^1_+$ with full support - for example
Gaussian measure in the half-line.
We will define inductively the measure $\nu$ on cone of distance matrices
$\cal R$ (indeed on ${\cal R}^0 $) by construction of its finite dimensional projection on
the cones ${\cal R}_n^0$ as follow.
The distribution of the entry $r_{1,2}$ is distribution $\gamma$.
So we have defined the measure on ${\cal R}_2^0$, denote it as $\nu_2$.
  Suppose we already define the joint distribution of the
entries $\{r_{i,j}\}_{i,j=1}^n$ which means that we define
a measure $\nu_n$ on ${\cal R}_n^0$.
 By  lemma 5 the open cone  ${\cal R}_{n+1}^0 $ is
homeomorphic to direct product of ${\cal R}_n^0$ and open half-space ${\bf R}^n_+$,
but we use only the structure of fibratoin:
projection ${\cal R}_n^0\stackrel{p_{n+1}}{\longleftarrow}{\cal R}_{n+1}^0$ -
gives us this structure with the fibers - interior of $A(r)$.
So in order to define a masure on  ${\cal R}_{n+1}^0 $ with  the given
measure $\nu_n$ on the base  ${\cal R}_n^0$ it is enough to define the conditional
measures on the fibers.
Choose each $r \in {\cal R}_n$ and represent the set of admissible vectors $A(r)$
as Minkowski sum: $$A(r)= M_r +\Delta_n,$$ (see 2.2) or as projection of the direct product
$\pi:M_r \times \Delta_n \to  M_r +\Delta_n=A(r)$. Consider product measure on
$M_r \times \Delta_n :
\gamma_r=m_r \times \gamma$ where $m_r$ is normalized Lebesgue measure on the compact
polytope $M_r$ and measure $\gamma$ was define before,
and its projection $\pi \gamma_r$ on $A(r)$. We define {\it conditional measure}
on $A(r)$ as $\pi \gamma_r$.

So we can define a measure on ${\cal R}_{n+1}$ denote it $\nu_{n+1}$
with given projection  $\nu_n$ on ${\cal R}_n^0$
and with the given set of conditional measures $\pi \gamma_r$ on the fibers $A(r)$.
By definition $\{\nu_n\}_{n=2}^\infty$
is a concordant systems of cylindric measures on ${\cal R}$ with respect to projections.
${\cal R}_n^0\stackrel{p_{n+1}}{\longleftarrow}{\cal R}_{n+1}^0$. Consequently
applying  Kolmogoroff's theorem we obtain the true measure $\nu$ on  ${\cal R}$.
The proof that for $\nu$-almost matrices  $r \in \cal R$ are universal repeats
the arguments of the proof of theorem 1, so $\nu \in \cal P$. This construction
indeed is pure combinatorial - we choose randomly a new points of metric space
and if we do this sufficiently ``random'' we obtain in completion the Urysohn space.
That is similar to the mentioned Erdos theorem about random graphs.
The detales will be done elsewhere.

We can varied the parameters of this construction, but of course it
 does not give in general a measure from the set $\cal D$ (matrix distribution)
because it could be not $S_{\bf N}$-invariant, and consequently not
from the set $\cal Q$. It means that these measures does not define
a measure on Urysohn space. What we did with this construction is
{\it a random countable every dense subset in Urysohn space}.
A very interesting problem is: how to obtain non degenerated measure on Urysohn space?

The natural but indirect method is to try to construct a measure on $\cal R$ which subtract
to all conditions of the theorem 5. We know a few things about geometry, topology and other
structures of Urysohn space $\cal U$ and about group of its
isometries in order to make a direct constructions.
Here is a concrete question about measures on it: is it possible to define a measure on
 $\cal U$ such that distance matrix of randomly and independently chosen $n$ points
has a given (for example Gaussian in the cone ${\cal R}_n$) distribution? Even the case $n=4$
is interesting. Another important problem - calculation of the matrix distributions for
concrete metric triples. Even finite dimensional distributions for compact Lie groups
with Haar measure or for some simple manifolds, as I know, never was found.

 \end{document}